\documentclass{amsart}
\usepackage{amssymb, amsmath, latexsym, amsthm}
\usepackage{fancyhdr}
\usepackage{graphicx}

\newtheoremstyle{dotless}{}{}{\itshape}{}{\bfseries}{}{ }{}
\theoremstyle{dotless}
\newtheorem{proposition}{Proposition}[section]

\newtheorem{theorem}{Theorem}[section]
\theoremstyle{definition}
\DeclareMathOperator{\supp}{supp}

\newcommand{\pd}{\partial}

\def\supp{\operatorname{supp}}

\newcommand{\al}{\alpha}

\newcommand{\la}{\lambda}

\newcommand{\om}{\omega}
\newcommand{\Om}{\Omega}

\newcommand{\cA}{\mathcal{A}}

\newcommand{\cD}{\mathcal D}
\newcommand{\cE}{\mathcal E}
\newcommand{\cF}{\mathcal F}

\newcommand{\cJ}{\mathcal J}

\newcommand{\cM}{\mathcal M}

\newcommand{\cR}{\mathcal R}
\newcommand{\cS}{\mathcal S}

\newcommand{\bI}{\mathbb I}

\newcommand{\bV}{\mathbb V}

\newtheorem{remark}[theorem]{Remark}
\newtheorem{defin}[theorem]{Definition}

\begin{document}

\title[Counterexamples for bi-parameter embedding]%
{%$H^1$ 
%Entropy bumps  and applications, the non-homogeneous setting
%Entropy ``bumping'' of the two weight Muckenhoupt $A_2$ condition
 %Bellman function sitting  on a tree
{Counterexamples for bi-parameter Carleson embedding}
%Calder\'on--Zygmund 
%operators
}
%\author[N.~Arcozzi]{Nicola Arcozzi}
%\address{Universit\`{a} di Bologna, Department of Mathematics, Piazza di Porta S. Donato, 40126 Bologna (BO)}
%\email{nicola.arcozzi@unibo.it}
%\thanks{Theorem 3.1 was obtained in the frameworks of the project 17-11-01064 by the Russian Science Foundation}
%\thanks{NA is partially supported by the grants INDAM-GNAMPA 2017 "Operatori e disuguaglianze integrali in spazi con simmetrie" and PRIN 2018 "Variet\`{a} reali e complesse: geometria, topologia e analisi armonica"}
%\author[I.~ Holmes]{Irina Holmes}
%\thanks{IH is partially supported by the NSF an NSF Postdoc under Award No.1606270}
%\address{Department of Mathematics, Michigan Sate University, East Lansing, MI. 48823}
\author[P. Mozolyako]{Pavel Mozolyako}
\thanks{PM is supported by the Russian Science Foundation grant 17-11-01064}
\address{Universit\`{a} di Bologna, Department of Mathematics, Piazza di Porta S. Donato, 40126 Bologna (BO)}
\email{pavel.mozolyako@unibo.it}
\author[G.~Psaromiligkos]{Georgios Psaromiligkos}
%\thanks{AV is partially supported by the NSF grant DMS-160065 and DMS 1900268 and by Alexander von Humboldt foundation}
\address{Department of Mathematics, Michigan Sate University, East Lansing, MI. 48823}
\email{psaromil@math.msu.edu \textrm{(G.\ Psaromiligkos)}}
\author[A.~Volberg]{Alexander Volberg}
\thanks{AV is partially supported by the NSF grant DMS-160065 and DMS 1900268 and by Alexander von Humboldt foundation}
\address{Department of Mathematics, Michigan Sate University, East Lansing, MI. 48823}
\email{volberg@math.msu.edu \textrm{(A.\ Volberg)}}
\makeatletter
\@namedef{subjclassname@2010}{
  \textup{2010} Mathematics Subject Classification}
\makeatother
% 42B	Harmonic analysis in several variables
% 42B20	Singular and oscillatory integrals (Calder?on-Zygmund, etc.)
% 42B35	Function spaces arising in harmonic analysis
% 47A	General theory of linear operators
% 47A30	Norms (inequalities, more than one norm, etc.)
%{30E20, 47B37, 47B40, 30D55.}
%
% 30D55	$H^p$-classes (1980-2009)
% 30E20	Integration, integrals of Cauchy type, integral representations of analytic functions
%
% 47B   	Special classes of linear operators
% 47B37	Operators on special spaces (weighted shifts, operators on sequence spaces, etc.)
% 47B40	Spectral operators, decomposable operators, well-bounded operators, etc.
\begin{abstract} 
We build here several counterexamples for two weight bi-parameter Carleson embedding theorem.
\end{abstract}
\maketitle

\section{Main definitions}

Let $T^2$ be a finite (but very deep) bi-tree. Bi-tree is the directed  graph  of all dyadic rectangles in the square $Q_0=[0,1]^2$. We assume that it terminates at small squares of size $2^{-N}\times 2^{-N}$, we call them generically by symbol $\om$.  The boundary $\pd(T^2)$ is this collection of $\om$'s. Often we identify $T^2$ with dyadic rectangles, whose family is called $\cD$. When we write $E\subset (\partial T)^2$ we mean any subset of $\om's$. It is convenient to think of $E$ as the union of ``$N$-coarse" dyadic rectangles.

Box condition
\begin{equation}\label{e:Boxcond}
\sum_{Q\in T^2, \,Q\subset R}\mu^2(Q)\alpha_Q \leq C_{\mu}\mu(E),\quad \textup{for any}\; R\in T^2.
\end{equation}
Carleson condition
\begin{equation}\label{e:Carlcond}
\sum_{Q\in T^2,\,Q\subset E}\mu^2(Q)\alpha_Q \leq C_{\mu}\mu(E),\quad \textup{for any}\; E\subset (\partial T)^2.
\end{equation}
Restricted Energy Condition
\begin{equation}\label{e:REC}
\sum_{Q\in \mathcal{D}}\mu^2(Q\cap E)\alpha_Q \leq C\mu(E),\quad \textup{for any}\; E\subset (\partial T)^2
\end{equation}
Embedding
\begin{equation}\label{e:imbed}
\sum_{Q\in\mathcal{D}}\left(\int_{Q}\varphi\,d\mu\right)^2\alpha_Q \leq C\int_{Q_0}\varphi^2\,d\mu\quad \textup{for any}\; \varphi\in L^2(Q_0,d\mu).
\end{equation}
Each of the next statement implies the previous one. We are interested when the opposite direction implications hold, and whether they hold in general.

Embedding is the boundedness of embedding operator $L^2(Q_0, \mu)\to L^2(T^2, \al)$ acting as follows:
\[
f\in L^2(Q_0, \mu)\to \Big\{\int_R f\, d\mu\Big\}_{R\in T^2}\in L^2(T^2, \al).
\]
 Restricted energy condition (REC)  is the  boundedness of embedding operator $L^2(Q_0, \mu)\to L^2(T^2, \al)$ on characteristic functions ${\bf 1}_E$ for all $E\subset\pd(T^2)$ uniformly.
 
 \section{Examples having box condition but not Carleson condition}
 \label{boxNotC}
 
 In \cite{Car} Carleson constructed the families $\cR$  of dyadic sub-rectangles of $Q=[0,1]^2$ having the following two  properties:
\begin{equation}
\label{boxCa}
\forall R_0 \in \cD, \quad \sum_{R\subset R_0, R\in \cR} m_2(R) \le C_0 m_2(R_0)\,,
\end{equation}
but for $U_\cR:= \cup_{R\in \cR} R$
\begin{equation}
\label{Ca}
 \sum_{ R\in \cR} m_2(R) \ge C_1 m_2(U_\cR)\,,
\end{equation}
where $C_1/C_0$ is as big as one wishes.

\medskip

 With {\it some} (rather wild) $\{\al_R\}_{R\in \cD}$ Carleson's counterexample is readily provide the example of $(\mu, \al)$ such that box condition is satisfied but Carleson condition fails. 
Indeed, let us put
$$
\al_R= \begin{cases} \frac{1}{m_2(R)}, \quad R\in \cR,
\\
0, \quad \text{otherwise}
\end{cases}
$$
Measure $\mu$ is just planar Lebesgue measure $m_2$. Fix any dyadic rectangle $R_0$, then box condition is satisfied:
$$
\sum_{R\subset R_0}\mu(R)^2 \al_R = \sum_{R\subset R_0, R\in \cR} m_2(R) \le m_2(R_0)=\mu(R_0)\,.
$$
But if $\Omega:= \cup_{R\in \cR} R$, then
$$
\sum_{R\subset \Omega}\mu(R)^2 \al_R =\sum_{R\subset \Omega, R\in \cR} m_2(R) =1\ge C\,m_2(\Omega)\,,
$$
where $C$ can be chosen as large as one wants. Hence, \eqref{Ca}  holds too with large constant.

\bigskip

The weight $\al=\{\al_R\}$. is rather wild here. But there is also an counterexample with $\al =1, 0$, see \cite{HPV}.

\section{Examples of having Carleson condition but not restricted energy condition REC.} 
Our aim here is to show that if we do not restrict ourselves to the constant weights as in \cite{AMPS18}, \cite{AHMV18b}, \cite{AMVZ19}, then the Carleson condition \eqref{e:Carlcond} is no longer sufficient for the embedding \eqref{e:imbed}. In fact even the Restricted Energy Condition \eqref{e:REC} is not necessarily implied by \eqref{e:Carlcond}. Namely we prove the following statement.
\begin{proposition}\label{p:71}
For any $\delta>0$ there exists number $N$, a weight $\alpha: T^2_N\mapsto \mathbb{R}_+$ and a measure $\mu $ on $(\partial T)^2$ such that $\mu$ satisfies Carleson condition \eqref{e:Carlcond} with the constant $C_{\mu} = \delta$,
\begin{equation}\label{e:771}
\sum_{Q\subset E}\mu^2(Q)\alpha_Q \leq \delta\mu(E),\quad \textup{for any }E\subset (\partial T)^2,
\end{equation}
but one can also find a set $F$ such that
\begin{equation}\label{e:772}
\sum_{Q\in\mathcal{D}}\mu^2(Q\cap F)\alpha_Q \geq \mu(F),
\end{equation}
hence the constant in \eqref{e:REC} is at least $1$.
\end{proposition}
We intend to give two examples of this kind. Both of them rely on the fact that one can basically consider this problem on a cut bi-tree by letting $\alpha$ to be either $1$ or $0$. 
This approach clearly does not work on the tree (see \cite{AHMV18a}[Theorem 1.1]), but the bi-tree has richer geometric structure. 

While globally (i.e. for $\alpha\equiv 1$) it looks at least somewhat similar to the tree (this similarity is implicit in the proof of \cite{AHMV18b}[Theorem 1.5]), one can remove some vertices (which is what essentially happens when we put $\alpha_Q := 0$) in such a way that the remaining part looks nothing like the full bi-tree (or a tree for that matter). In particular this allows us to create a significant  between the amount of ``available" rectangles that lie inside $E$ or just intersect $E$ for a certain choice of the weight $\alpha$ and the set $E$ (this corresponds to the difference between Carleson and REC conditions).\par
The first example is quite simple and is inspired by the counterexample for $L^2$-boundedness of the biparameter maximal function. The weight $\alpha$ in this case cuts most of the bi-tree, and the resulting set differs greatly from the original graph.\\ The second example is somewhat more involved, on the other hand the weight there leaves a much bigger portion of the bi-tree, actually it has a certain monotonicity property: $\alpha_R \geq \alpha_Q$ for $R\supset Q$. The structure of the ``available" set is more rich in this case (it looks more like $\mathbb{Z}^2$ in a sense), nevertheless there is not enough rectangles to have the Carleson-REC equivalence.

\subsection{A simple example of having Carleson condition but not restricted energy condition}
Let $N\in\mathbb{N}$ be some large number (to be specified later), and let $T^2 = T_N^2$ be a bi-tree of depth $N$. We use the dyadic rectangle representation of $T^2$.\par
Let $\om$ be our $[0,2^{-N}]^2$ left lower corner. Given $R = [a,b]\times[c,d]\in\mathcal{D}$ let $R^{++} := \left[\frac{a+b}{2},b\right]\times\left[\frac{c+d}{2},d\right]$ be the upper right quadrant of $R$. Consider
$Q_1= [0,1]\times [0, 2^{-N+1}]$ and its $Q_1^{++}$, and
$Q_2= [0,2^{-1}]\times [0, 2^{-N+2}]$ and its $Q_2^{++}$,  
$Q_3= [0,2^{-2}]\times [0, 2^{-N+3}]$ and its $Q_3^{++}$,  et cetera\dots . In total, $N$ of them.

Put measure $\mu$ to have mass $\tau_0:=1/\sqrt{N}$ on $\om$, and uniformly distribute mass $\tau_i$ on $Q_i^{++}$.

Now $\al_R$ is always zero except when $R= \om, Q_1, Q_2, \dots$. For those $\al=1$, so we have $N+1$ alphas equal to $1$.

Now choose set $E=\om$. When we calculate $\cE[\mu|E]$ we sum up
$$
\tau_0^2 +\sum_{i=1}^N \mu(\om\cap Q_i)^2 =( N+1)\tau_0^2 =\frac{N+1}{N}\ge \sqrt{N} \frac1{\sqrt{N}} =\sqrt{N}\mu(E)\,.
$$

 So REC condition has a big constant.
 
 \bigskip

%On the other hand if $\Om=\cup R_j$ for some $R_j$, then it is useless to consider $R_j$ that are not $\om, Q_1, Q_2,...$ since their contribution into $\cE_\Om[\mu]$ is zero as 
Since
$$
\cE_\Om[\mu] =\sum_{R\subset \Om,\, \al_R\neq 0} \mu(R)^2\,,
$$
then, denoting $Q_0:=\om$, we have
$$
\cE_\Om[\mu] = \sum_{j: Q_j\subset \Om} \mu(Q_j)^2 =:\sum_{j\in J( \Om)} \mu(Q_j)^2\,.
$$

%So, for example if $\Om$ as a set does not include all $\om, Q_1, Q_2,...$, then $\cE_\Om[\mu] =0$.

Let $\tau_0=1/\sqrt{N} \le \frac14=\tau_1=\tau_2=...=\tau_N$.
Then  
$$
\cE_\Om[\mu] =\sum_{j\in J(\Om)} (\tau_0+\tau_j)^2 \le 4\sum_{j\in J(\Om)} \tau_j^2\,.
$$
And this is $\le \sum_{j\in J(\Om)} \tau_j \le \mu(\Om)$.
So Carleson condition holds with constant $1$.

\subsection{The lack of maximal principle matters}
\label{maxP}

All measures and dyadic rectangles  below will be $N$-coarse.

In this section we build another example when Carleson condition holds, but restricted energy condition fails. 
But the example is more complicated (and more deep) than the previous one. In it the weight $\al$  again has values either $1$ or $0$, but the support $S$ of $\al$ is an {\it up-set}, that is, it contains every ancestor of every rectangle in $S$. 

The example is based on the fact that potentials on bi-tree may not satisfy maximal principle. So we start with constructing $N$-coarse $\mu$ such that given a small $\delta>0$
\begin{equation}
\label{le1}
\bV^\mu \lesssim \delta\quad \text{on}\,\, \supp\mu,
\end{equation}
but  with an absolute strictly positive $c$
\begin{equation}
\label{logN}
\max \bV^\mu \ge \bV^\mu(\om_0) \ge c\, \delta\log N\,,
\end{equation}
where $\om_0:=[0, 2^{-N}]\times [0, 2^{-N}]$.

We define a collection of rectangles
\begin{equation}\label{e:773}
Q_j :=  [0,2^{-2^j}] \times  [0,2^{-2^{-j}N}] ,\quad j=1\dots M\approx \log N,
\end{equation}
and we let 
\begin{equation}\label{Qs}
\begin{split}
&Q_j^{++} :=  [2^{-2^j-1},2^{-2^j}]\times  [2^{-2^{-j}N-1},2^{-2^{-j}N}] \\
&Q_j^- := [0,2^{-2^j-1}] \times  [0,2^{-2^{-j}N}] ,\quad j=1\dots M\\
&Q_j^t := Q_j\setminus Q_j^-\\
&Q_j^r:= [2^{-2^j-1},2^{-2^j}]\times  [0,2^{-2^{-j}N}]\\
& Q_j^{--} := Q_j^-\setminus Q_j^r
\end{split}
\end{equation}
to be their upper right quadrants, lower halves, top halves, right halves, and lower quadrant respectively.
Now we put 
\begin{equation}\label{e:774}
\begin{split}
&\mathcal{R} := \{R:\; Q_j\subset R\; \textup{for some }j=1\dots M\}\\
&\alpha_{Q} := \chi_{\mathcal{R}}(Q)\\
%&\mu(\omega_0) := \frac{1}{MN}\\
&\mu(\omega) := \frac{\delta}{N}\sum_{j=1}^M \frac{1}{|Q_j^{++}|}\chi_{Q_j^{++}}(\omega),\\
%&F := \bigcup_{j=1}^M Q^-_j,
& P_j= (2^{-2^j}, 2^{-2^{-j}N})\,.
\end{split}
\end{equation}
here $|Q|$ denotes the total amount of points $\omega \in (\partial T)^2\cap Q$, i.e. the amount of the smallest possible rectangles (of the size $2^{-2N}$) in $Q$.\par
%We claim that the measure, 
%weight and set defined above satisfy the claim of Proposition.
Observe that on $Q_j$ the measure is basically a uniform distribution of the mass $\frac{\delta}{N}$ over the upper right quarter $Q_j^{++}$ of the rectangle $Q_j$ (and these quadrants are disjoint).
\par

To prove \eqref{le1} we fix $\om\in Q_j^{++}$ and split $\bV^\mu(\om)= \bV^\mu_{Q_j^{++}}(\om) +\mu(Q_j^t)+ \mu(Q_j^r)+ \bV^\mu(Q_j^{++})$, where the first term sums up $\mu(Q)$ for $Q$ between $\om$ and $Q_j^{++}$. This term obviously  satisfies $\bV^\mu_{Q_j^{++}}(\om) \lesssim \frac{\delta}{N}$. Trivially $\mu(Q_j^t)+ \mu(Q_j^r)\le  \frac{2\delta}{N} $. The non-trivial part is the estimate
\begin{equation}
\label{VQ}
\bV^\mu(Q_j^{++}) \lesssim \delta\,.
\end{equation}
To prove \eqref{VQ}, consider the sub-interval of interval $[1, n]$ of integers. We assume that $j\in [m, m+k]$.  We call by $C^{[m, m+k]}_j$ the family of dyadic rectangles containing $Q_j^{++}$ along with all $Q_i^{++}$, $i\in [m, m+k]$ (and none of the others). Notice that $C^{[m, m+k]}_j$ are not disjoint families, but this will be no problem for us as we  wish to
estimate $\bV^\mu(Q_j^{++})$ from above. 

Notice that, for example, $C^{[m, m+1]}_j$ are exactly the dyadic rectangles containing point $P_j$. It is easy to calculate that the number of such rectangles
is 
\[
(2^j+1)\cdot (2^{-j}N +1)\lesssim N\,.
\]

Analogously, dyadic rectangles in family $C^{[m, m+k]}_j$ have to contain  points $P_m, P_{m+k}$. Therefore,  each of such rectangles contains  point $(2^{-2^m}, 2^{-2^{-m-k}N})$.
The number of such rectangles is obviously at most $\lesssim 2^{-k}N$. The number of classes $C^{[m, m+k]}_j$ is at most $k+1$.

Therefore, $\bV^\mu((Q_j^{++})$ involves at most $(k+1)2^{-k}N$ times the measure in the amount $k\cdot \frac{\delta}{N}$. Hence
\[
\bV^\mu((Q_j^{++}) \le \sum_{k=1}^n k(k+1)2^{-k}N \cdot \frac{\delta}{N}\,,
\]
and \eqref{VQ} is proved. Inequality \eqref{le1} is also proved.

\medskip

We already denoted
\[
\om_0:=[0, 2^{-N}]\times [0, 2^{-N}]\,,
\]
 calculate now $\bV^\mu(\om_0)$. In fact, we will estimate it from below. The fact that $C^{[m, m+k]}_j$ are not disjoint may represent the problem now because  we wish estimate $\bV^\mu(\om_0)$ from below.
 
 To be more careful for every $j$ we denote now by $c_j$ the family of dyadic rectangles containing the point $P_j$ but not containing any other point
$P_i, i\neq j$. Rectangles in $c_j$ contain $Q_j^{++}$ but do not contain any of $Q_i^{++}$, $i\neq j$. There are 
$(2^{j-1} - 2^j-1)\cdot (2^{-j+1}N-2^{-j}N -1)$, $j=2, \dots, M-2$. This is at least $\frac18 N$.

But now families $c_j$ are disjoint, and rectangles of class $c_j$ contribute at least $\frac18 N\cdot \frac{\delta}{N}$ into the sum that defines $\bV^\mu(\om_0)$. W have $m_4$ such classes $c_j$, as $j=2, \dots, M-2$.
Hence,
\begin{equation}
\label{deltaM}
\bV^\mu(\om_0) \ge \frac18 N\cdot \frac{\delta}{N}\cdot (M-4)\ge\frac19 \delta M\,.
\end{equation}
Choose  $\delta$   to be a small absolute number $\delta_0$. Then we will have (see \eqref{le1})
\[
\bV^\mu \le 1, \quad\text{on}\,\, \supp\mu\,.
\]
But \eqref{deltaM}  proves also \eqref{logN} as $M\asymp \log N$.

\medskip

\begin{remark}
\label{Mx}
Notice that in this example $\bV^\mu\le 1$ on $\supp\mu$, and
\begin{equation}
\label{exp}
\text{cap} \{\om: \bV^\mu\ge \la\} \le c e^{-2\la}\,.
\end{equation}
Here capacity is the bi-tree capacity defined e. g. in \cite{AMPS18}. So there is no maximal principle for the bi-tree potential, but the set, where the maximal principle breaks down,  has small capacity.
\end{remark}

\bigskip

Now we construct the example of $\nu$ and $\al$ with $\al=1$ on an up-set (and zero otherwise), and such that Carleson condition is satisfied
but REC (restricted energy condition) is not satisfied. We use the same measure $\mu$ we have just constructed, and we put
\[
\nu := \mu +\nu|\om_0,
\]
where $\nu|\om_0$ is the uniformly distributed over $\om_0$ measure of total mass $\frac1{MN}$. Weight $\al$ is chosen as in \eqref{e:774}.

\bigskip

%%%%%%%%%%%%%%%%%%%%%%%%
\noindent{\bf Warning.} The meaning of $\bV$ changes from now on. Before $\bV^{\cdot}= \bI \bI^*(\cdot)$. Everywhere below,
\[
\bV^\nu:= \bI[\al \bI^*(\nu)]\,.
\]
%%%%%%%%%%%%%%%%%%%%%

\bigskip

Let us first check that REC constant  is bad. We choose $F=\cup_j Q_j^-$. Then $\nu_F:=\nu|F = \nu|\om_0$. On the hand, and this is the main feature,
%%%%%%%%%%%%%%%%%%%%%%
\begin{equation}
\label{lies_in_M}
\om_0 \,\,\text{lies in}\,\, M\,\, \text{rectangles} \,\, Q_j\,.
\end{equation}
%%%%%%%%%%%%%%%%%%%%%%%%%
Hence,  
%%%%%%%%%%%%%%%%%%%%%%%
\begin{equation}
\label{number_of_rectangles}
\text{there are} \,\, \ge c NM\,\, \text{dyadic rectangles}\,\, R\,\,\text{ such that}\,\, \al_R=1\,\, \text{ and}\,\, \om_0\subset R\,. 
\end{equation}
%%%%%%%%%%%%%%%%%%%%%%%
 In fact, consider dyadic rectangles  in $\cup_j c_j$, where families $c_j$ were built above. For each
$R\in \cup c_j$ we have $\al_R=1$, see \eqref{e:774}. And there are at least $\frac18 NM$ of them.  We conclude
%%%%%%%%%%%%%%%%%%%
\begin{equation}
\label{bV_at_om0}
\bV^{\nu_F}(\om_0) \ge \frac18 MN\cdot \nu(\om_0)\,. 
\end{equation}
%%%%%%%%%%%%%%%%%%%%%%%%%
Therefore, 
\[
\int \bV^{\nu|F}\, d \nu|F \ge  \nu(\om_0)^2\cdot \frac18 NM=\frac18 \frac1{MN} \ge c_0\,\nu(\om_0)\,.
\]
This means that constant of REC is at least absolute constant $c_0$. Let us show that the Carleson constant is $\lesssim c\cdot \delta$. But $c_0$ has nothing to do with $\delta$ that can be chosen as small as we wish. 
\begin{remark}
\label{inM}
We do not need the following claim  now, we will need it only later, but notice that in a fashion completely similar to the one that just proved \eqref{bV_at_om0}, one can also prove
\begin{equation}
\label{bVmu_at_om0}
\bV^{\mu}(\om_0) \ge \frac18 MN\cdot \frac{\delta}{N} \ge c\, \delta M\,. 
\end{equation}
Moreover, we already proved it in \eqref{deltaM}. This holds because $\om_0$ is contained in exactly $M$ rectangles $Q_j^{++}$.
\end{remark}

\bigskip

\begin{defin}
Dyadic rectangles whose left lower corner is $(0,0)$ will be called {\it hooked} rectangles.
\end{defin}

To check the Carleson condition with small constant we fix any finite family $\cA$ of dyadic   rectangles, and let
\[
A=\cup_{R\in \cA} R\,.
\]

We are interested in subfamily $\cA'$ of $R$ such that $\al_R=1$. Other elements of $\cA$ do not give any contribution to $\cE_A[\nu]$ as $ \mu(Q)^2 \al_Q=\mu(Q)^2 \cdot  0=0$ for any $Q\subset Q', \, Q'\in \cA\setminus \cA',$  as the support of $\al$ is an up-set.

All rectangles from $\cA'$ are hooked rectangles. As we noticed, we can think that $\cA'=\cA$.  In other words, without the loss of generality, we can think that $\cA$ consists only of hooked rectangles. Any hooked rectangle generates  a closed interval $\cJ$ in the segment $[1, n]$ of integers: interval $\cJ$ consists of $j$, $1\le j\le n$, such that point $P_j$ lies in this hooked rectangle.  This is the same as to say that $Q_j$, $j\in \cJ$, is a subset of this hooked rectangle. 

So family $\cA$  generates the family of  closed intervals in the segment $[1, n]$ of integers. Let us call $\cJ_A$ this family of intervals in the segment $[1, n]$ of integers. Intervals of family $\cJ_A$ can be not disjoint. But we can do the following, if intervals intersect, or even if these closed intervals are adjacent, we unite them to a new interval. The new system (of disjoint and not even adjacent)  closed intervals corresponds to another initial system $\tilde \cA$, and we can think that $\tilde \cA$ consists of hooked rectangles. 
We call a system of hooked rectangles  {\it a clean system}  if it gives rise to not adjacent disjoint family of closed intervals inside the set $[1,n]$ of integers
The relationship between rectangles in $\cA$ and $\tilde\cA$ is the following: each rectangle of $\tilde\cA$ is a common ancestor of a group of rectangles in $\cA$.  

A very important geometric property of $\tilde \cA$ is the following. Let $Q\in \tilde\cA$ and let $R^1, \dots, R^s$ be all rectangles from $\cA$ such that $R^i\subset Q$, $i=1,\dots, s$. Then 
\begin{equation}
\label{tildeA}
\nu(Q\setminus \cup_{i=1}^m R^i) =0\,.
\end{equation}

In particular, \eqref{tildeA} implies
\begin{equation}
\label{tildeA1}
\nu(A)=\nu(\tilde A) \,.
\end{equation}

When checking the Carleson condition
\begin{equation}
\label{carlA}
\cE_A[\nu] \lesssim \delta \nu(A),
\end{equation}
we can always think about $\cA$ being replaced by $\tilde \cA$ and $A$ being replaced by $\tilde A$  because in \eqref{carlA} the RHS stays the same, but the LHS can jump up
with passage from $A$ to $\tilde A$. Therefore, checking \eqref {carlA} for {\it clean systems} of rectangles is the same as to check it for all systems of rectangles. From now on $\cA$ is supposed to be clean.

To prove \eqref{carlA} is the same as to prove (since $\nu_A(Q) = (\mu_A(Q) +\nu_{\om_0}(Q)$)
\begin{equation}
\label{carlA1}
 \sum_{Q\subset A} \nu_A(Q)^2 \al_Q  + \sum_{Q\subset A} \nu_{\om_0}(Q)^2 \al_Q  \lesssim \delta \nu(A),
\end{equation}

The first sum is bounded by $\int \bV^{\mu|A}\, d\mu|A$. But by \eqref{le1} (which  follows from \eqref{VQ}) we have
\[
\sum_{Q\subset A} \nu_A(Q)^2 \al_Q = \int \bV^{\mu|A}\, d\mu|A \lesssim \delta \|\mu_A\|= \delta \mu(A) \le \delta \nu(A)\,,
\]
and this means that the part of \eqref{carlA1} is proved.

To estimate $\sum_{Q\subset A} \nu_{\om_0}(Q)^2 \al_Q= \frac1{(MN)^2} \sharp\{R: \al_R=1, R\subset A\} $  we take one interval $\cJ_k$ from the family generated by the clean $\cA$ in $[1,n]$,  we denote
\begin{equation}
\label{mk}
m_k:=\sharp\cJ_k\,,
\end{equation}
and we estimate how many dyadic rectangles  $R$ contain one of $Q_j$, $j\in \cJ_k$. We even do not care now whether $R$ is a subset of $A$ or not.  The number of such rectangles in at most $m_k\cdot N$. On the other hand,
\[
\nu(A)\ge \sum_k m_k \cdot \frac{\delta}{N} + \frac1{MN} \sum_k m_k \ge \sum_k m_k \cdot \frac{\delta}{N} \,.
\]
So to prove the estimate for the second sum in \eqref{carlA1}, we need to see that
\[
 \frac1{(MN)^2}\sum_k m_k\cdot N \lesssim \frac1M\sum_k m_k \cdot \frac{\delta}{N},
 \]
 which is obviously true if we choose $\delta \ge \frac1{M}=\frac1{\log N}$. So \eqref{carlA1} is proved.

 %%%%%%%%%%%%%%%%%%%%%%%%%%%%%%%%%%%%%%%%
 %%%%%%%%%%%%%%%%%%%%%%%%%%%%%%%%%%%%%%%
% \begin{remark}
% \label{product}
% The system $(\nu, \al)$ constructed above gives the example of measure and $\al=1, 0$ with 
%support on an up-set (so on a nice subgraph of a bi-tree) such that Carleson condition holds, 
%but REC condition does not hold. Notice that the support of $\al$ consists of 
%$m\approx \log N$  sets $S_1, \dots, S_M$ in bi-tree, where each set $S_i$ can be given by the equation
% \begin{equation}
 %\label{Si}
% S_i=\{ Q=I\times J: \tau_i(I)\cdot \eta_i(J) =1\}\,,
% \end{equation}
 %where $\tau, \eta$ are characteristic functions of certain up-sets in a simple tree $T$.  So
% \[
% \al=\sum_{i=1}^M \tau_i\otimes \eta_i\,.
% \]
% The discrepancy between Carleson constant and REC constant is at least of the order $\frac1{M}$.
%But its {at least}, nobody said that it is equivalent to $1/M$, where $M$ is the least number in display formula above.  
%To understand better the relationship between the least number in display formula above and the 
%REC-Carleson discrepancy one can ask what happens if $M$ is small, like $M=1$ 
%(this means that  the weight $\al$ has a product structure), or $M\asymp 1$? 
%It turns out that in this case $REC$ and Carleson conditions are equivalent, see \cite{AMVZ19}.
% \end{remark}
 %%%%%%%%%%%%%%%%%%%%%%%%%%%%%%%%%%%%%%%%
 %%%%%%%%%%%%%%%%%%%%%%%%%%%%%%%%%%%%%%%

\section{Restricted energy condition holds but no embedding}
\label{RECnoEmb}

In this section we emulate the previous construction, we start with $\{Q_j\}$ and measure $\mu$ but instead of adding $\om_0$ we will add a more sophisticated piece of measure.

Let us start with  recalling the system $\{Q_j\}, j=1, \dots, M$ and measure $\mu$ from the previous section. We continue with denoting 
\[
Q_{0, j}:= Q_j,\quad \mu_0:=\mu\,\, \text{from the previous section}\,.
\]
Next we continue with defining a a sequence of collections $\mathcal{Q}_k,\; k=0\dots K,$ of dyadic rectangles as follows
\begin{equation}\label{e:874}
\begin{split}
\mathcal{Q}_k := \left\{Q_{k,j} = \bigcap_{i=j}^{j+2^k}Q_{0,i},\; j=1\dots M-2^k\right\},\; k=1\dots K.
\end{split}
\end{equation}
In other words, $\mathcal{Q}_0$ is the basic collection of rectangles, and 
%%%%%%%%%%%%%%%%%%%%%%
\begin{equation}
\label{lies_in_2k}
\mathcal{Q}_{k}\,\text{ consists of the intersections of} \,\,2^k\,\text{ consecutive elements of}\,\mathcal{Q}_0.
\end{equation}
%%%%%%%%%%%%%%%%%%%%%%%%%%%%%%%%%%%%%%%%%%%%%%%%%%%%%%
 The total amount of rectangles in $\mathcal{Q}_k$ is denoted by $M_k = M-2^{k}+1$.

We also denoted by $\cR$ the collection of rectangles lying above $\mathcal{Q}_0$
\begin{equation}\label{e:874.4}
\cR := \{R:\; Q_{0,j}\subset R\;\textup{for some}\; 1\leq j\leq M\},
\end{equation}
and we let
\[
S_k := \bigcup_{Q\in \mathcal{Q}_k}Q.
\]
The weight $\alpha$ was defined as follows:
\begin{equation}\label{e:874.5}
\begin{split}
&\alpha_{Q} := 1,\quad \textup{if}\; Q\in \cR \\
&\alpha_Q := 0 \quad \textup{otherwise}.
\end{split}
\end{equation}

Now we construct the measure $\mu$, whose main part will be already constructed $\mu_0$. Let 
\[
Q_{k,j}^{++} := \left[2^{-2^{j+2^k}-1},2^{-2^{j+2^k}}\right]\times \left[2^{-2^{j}N-1},2^{-2^{j}N}\right]
\]
be the upper right quadrant of $Q_{k,j}$. For every $k = 0\dots K$ we distribute the mass $2^{-2k}\frac{M_k\delta}{N}$ over the rectangles $Q_{k,j}^{++}$. Namely, for every $j=1\dots M_k$ we attach a mass $2^{-2k}\frac{\delta}{N}$ to the rectangle $Q_{k,j}$ that is uniformly distributed over the quadrant $Q_{k,j}^{++}$. We note that all these quadrants $Q_{k,j}^{++}$ are disjoint.\par
Measure $\mu_0$ is the ``main" part of $\mu$, in the sense that $\mu_0$ is generated by the masses on $\mathcal{Q}_0$,
\[
\mu_0(\omega) := \frac{\delta}{N}\sum_{j=1}^M \frac{1}{|Q_{0,j}^{++}|}\chi_{Q_{0,j}^{++}}(\omega),\quad \omega \in (\partial T)^2,
\]
and let $\mu_k$ be the corresponding mass on $\mathcal{Q}_k$
\[
\mu_k(\omega) := \frac{2^{-2k}\delta}{N}\sum_{j=1}^{M_k} \frac{1}{|Q_{k,j}^{++}|}\chi_{Q_{0,j}^{++}}(\omega),\quad \omega \in (\partial T)^2,
\]
so that 
\[
\mu = \mu_0 + \sum_{k=1}^{K}\mu_k.
\]
Finally we define the function $f$, and we do it in such a way that it is 'congruent' with the distribution of $\mu_0$ over $\Omega$, namely we let
\[
f(R) := \mu_0(R)\cdot \alpha_R.
\]

\subsection{Main idea}
\label{MI}

Notice that 
\[
\bV^{\mu_0} = \bI f=\bI[\al \bI^*\mu_0], \, \,\,\int \bV^{\mu_0}\, d\mu_0 =\sum\bI^*\mu_0\cdot \bI^*\mu_0 \cdot\al= \sum_{T^2} f^2\cdot \al\,.
\]
To prove that embedding has a bad constant, it is sufficient to show that the dual inequality has a bad constant:
\[
\int (\bI f)^2 \, d\mu >> \sum f^2\cdot \al,
\]
which becomes 
\begin{equation}
\label{dual}
\int (\bV^{\mu_0})^2\, d\mu>> \int \bV^{\mu_0} \, d\mu_0\,.
\end{equation}

Let us look at Remark \ref{inM}, at 
\eqref{lies_in_M}, \eqref{number_of_rectangles}, 
\eqref{bVmu_at_om0} and compare \eqref{lies_in_M} with  \eqref{lies_in_2k}.
The conclusion is: since every $Q_{k,j}$ lies in $2^k$ of $Q_{0,j}$ (number $2^k$ replaces $M$ in \eqref{lies_in_M}, \eqref{number_of_rectangles}, \eqref{bVmu_at_om0}),  then
\begin{equation}
\label{bVmu0_on_muk}
\bV^{\mu_0} \ge c\,2^k N \cdot \frac{\delta}{N}=c  \delta 2^k\,\,\text{on each}\,\, Q_{k, j}\,. 
\end{equation}

We already saw that $\bV^{\mu_0}\lesssim \delta$ on $\mu_0$, so
\begin{equation}
\label{RHS}
\int \bV^{\mu_0} \, d\mu_0 \lesssim \delta^2 \frac{M}{N}\,.
\end{equation}

Now, using \eqref{bVmu0_on_muk}
we get
\begin{equation}
\label{LHS}
\int (\bV^{\mu_0})^2\,d\mu =\sum_{k=1}^K \int (\bV^{\mu_0})^2\,d\mu_k \ge c^2 \delta^2 \sum_{k=1}^K  2^{2k} \|\mu_k\| = c^2 \delta^3 \frac{M\log M}{N}
\end{equation}
 
 For example, let
 \begin{equation}\label{choice}
 \delta= \frac1{\log M}\,.
 \end{equation}
 Then the constant of  embedding is $\approx 1$.
 
 \subsection{REC condition holds with a small constant} 
 \label{RECsmall}
 
 Let $\cA$ be a collection of (hooked) rectangles, $A=\cap_{R\in \cA} R$. Let $\nu_k:=\mu_k|A$, $k=1, \dots, K$, $\nu:= \mu |A=\sum \nu_k$.
 We need to  prove
 \begin{equation}
 \label{smREC}
 \cE_A[\nu]\lesssim \delta \|\nu\|\,.
 \end{equation}
 Let $n>k$, we wish to estimate $\bV^{\mu_n}(Q_{k, j})$. This is a certain sum over a system $\cS$ of rectangles  of the form
 \[
 \sum_{R\in \cS} \mu_n(R),
 \]
 where dyadic rectangles are a) contain $Q_{k, j}$, b) $\al_R=1$. Notice that this system depends on $Q_{k, j}$ but totally independent of $n$. So if we manage to estimate $\mu_n(R)$ via $\mu_0(R)$, then we compare $\bV^{\mu_n}(Q_{k, j})$ to $\bV^{\mu_0}(Q_{k, j})$.
 
 But let the number of $Q_{s, j}$ in $R$ be denoted by $m^s_R$. Then it is very easy to see that
 \[
 m^n_R \le m^0_R +2^n+1\,.
 \]
 Then 
 \[
 2^{-2n}m^n_R \le 2^{-n}(m^0_R +1+2^{-n}) \le 3\cdot  2^{-n} m^0_R\,.
 \] 
 Then 
  \[
\bV^{\mu_n}(Q_{k, j}) =\sum_{R\in \cS} \mu_n(R) \le 3\cdot  2^{-n} \sum_{R\in \cS} \mu_0(R) = 3\cdot 2^{-n} \bV^{\mu_0}(Q_{k,j})\lesssim \delta 2^{k-n} \,.
 \]
 Therefore,
\[
\sum_{n\ge k} \bV^{\nu_n} \, d\nu_k\lesssim \delta \sum_{n\ge k} 2^{k-n} \|\nu_k\|=2\delta \|\nu_k\|\,.
\]
And so
\[
\sum_k\sum_{n\ge k} \bV^{\nu_n} \, d\nu_k\lesssim \delta \sum_k \|\nu_k\|=\delta\|\nu\|\,.
\]
Inequality \eqref{smREC} is proved.

%%%%%%%%%%%%%%%%%%%%%%%%%%%%%%%%%%%%
%%%%%%%%%%%%%%%%%%%%%%%%%%%%%%%%%%%%%
%\begin{remark}
%\label{logMdiscrepancy}
%We have
 %\[
% \al=\sum_{i=1}^M \tau_i\otimes \eta_i\,.
% \]
% The discrepancy between REC constant and embedding constant is at least of the order $\frac1{\log M}$.
%\end{remark}
%%%%%%%%%%%%%%%%%%%%%%%%%
%%%%%%%%%%%%%%%%%%%%%%%%%%%%%

\section{Strong dyadic maximal function and counterexamples}
\label{counter}

%We start to consider two weight bi-parameter Carleson embedding. We  have an arbitrary finite positive measure (weight) on $Q_0:=[0,1]^2$ and a sequence of non-negative %numbers $\al_Q$ enumerated by the dyadic sub-rectangles of $Q_0$ (=vertices of a bi-tree).

\noindent{\bf Definition.}
Let $\cS$ be a family of dyadic sub-rectangles of $Q_0$ (may be $\cS=\cD$, the family of all dyadic sub-rectangles). We call the sequence of $\{\beta_Q\}_{Q\in \cS}$ Carleson if
\begin{equation}
\label{bC}
\forall \cS'\subset \cS, \quad \sum_{Q\in \cS'}\beta_Q \mu(Q) \le C\mu(\cup_{Q\in \cS'} Q)\,.
\end{equation}
The best $C$ is called the Carleson norm of the sequence.

\noindent{\bf Definition.} Abusing the language we say that the weight $\al:=\{\al_Q\}_{Q\in \cS}$ satisfies Carleson condition if the sequence $\beta_Q:= \al_Q \cdot \mu(Q)$ is a Carleson sequence:
\begin{equation}
\label{alC}
\forall \cS'\subset \cS, \quad \sum_{Q\in \cS'}\al_Q \mu(Q)^2 \le C\mu(\cup_{Q\in \cS'} Q)\,.
\end{equation}

We already know, see \cite{AHMV18b} e.g., that bi-parameter Carleson embedding
\[
\sum_{R\in \cD} \big(\int_R \psi\, d\mu\big)^2 \al_R =\sum_{R\in \cD} [\bI^*(\psi \mu)]^2 \al_R  \le C' \int_{Q_0} \psi^2\, d\mu
\]
is equivalent to the Carleson condition above if we have $\al_R=1,\hspace{0.1cm} \forall R \in \cD$ .

Understanding the general two-weight bi-parameter situation (that is $\al \not\equiv 1$) seems to be super hard, as the examples above show. Notice that in two-weight  one-parameter situation the answer is known, see, e.g. \cite{NTV99}. And the answer is given in terms of Carleson condition.  However, in bi-parameter situation this is far from being so as the following theorem shows. First we give

\noindent{\bf Definition.} A finite positive measure $\mu$ on $Q_0$ is called a ``bad'' measure if there exists
weight $\al=\{\al_Q\}_{Q\in \cD}$ that satisfies the Carleson condition but such that the embedding
\begin{equation}
\label{alEmb}
\sum_{R\in \cD} \Big(\int_R \psi\, d\mu\Big)^2\al_R \le C' \int_{Q_0} \psi^2\, d\mu
\end{equation}
does not hold.

The strong maximal function with respect to $\mu$  is
\[
\cM_\mu \psi (x) =\sup_{\substack{R\in \cD \\ R\ni x}} \frac{1}{\mu(R)}\int_R |\psi|\, d\mu\,,
\]
where $0/0=0$. The supremum is taken over all dyadic sub-rectangles of $Q_0$.

\begin{theorem}
\label{Dor}
Let $\mu$ be atom free. Then the measure $\mu$  is bad if and only if $\cM_\mu$ is not a bounded operator in $L^2(\mu)$.
\end{theorem}

We need some preparation, which is written down in \cite{Verb} and \cite{TH}, but we repeat it for the convenience of the reader. 

If weight $\al$ satisfies Carleson condition with constant $1$, then $\beta_Q=\al_Q \mu(Q)$ is a Carleson sequence with Carleson constant $1$ and this means that for any family $\cD'\subset \cD$
\begin{equation}
\label{alCarl1}
\sum_{Q\in \cD'} \al_Q \mu(Q)^2 \le \mu(\bigcup_{Q\in\cD'} Q)\,.
\end{equation}

Using the idea of Igor Verbitsky \cite{Verb} this can be stated in terms of discrete Triebel-Lizorkin space $f^{1, \infty}(\mu)$:
\[
\{\al_Q\mu(Q)\}_{Q\in \cD}\in f^{1, \infty}(\mu),
\]
and has norm at most $1$. But $( f^{1, \infty}(\mu))^* =  f^{ \infty, 1}(\mu)$, and this is the space of coefficients $\{\lambda_Q\}_{Q\in \cD}$ given by the norm
\[
\|\{\lambda_Q\}_{Q\in \cD}\|_{ f^{ \infty, 1}(\mu)} := \int\sup_{Q\in\cD}\Big(\lambda_Q\chi_Q(x)\Big)\, d\mu\,.
\]
Thus, by duality \eqref{alCarl1} is equivalent to
\begin{equation}
\label{alCarlDual}
\forall \{\lambda_Q\}_{Q\in \cD}, \quad\sum_{Q\in \cD} \al_Q \mu(Q)^2\cdot \lambda_Q \le \int \sup_{Q\in\cD}\Big(\lambda_Q\chi_Q(x)\Big)\, d\mu\,.
\end{equation}
 Without loss of generality we think that all $ \al_Q \mu(Q)>0$ (otherwise they are not in the LHS and we forget them). Then we can rewrite \eqref{alCarlDual} as
 
\begin{equation}
\label{alCarlDual2}
\forall \{b_Q\}_{Q\in \cD}, \quad \sum_{Q\in \cD} b_Q \le \int \sup_{Q\in\cD}\Big(b_Q\frac{\chi_Q(x)}{\al_Q \mu(Q)^2}\Big)\, d\mu\,.
\end{equation}

In \cite{Verb} and in \cite{TH} now the theorem of L. Dor \cite{Dor} is used. Here it is
\begin{theorem}
\label{thDor}
Let $\mu$ be atom free. Let $g_i$ be a fixed sequence of measurable functions and suppose for all positive numbers $b_i$ that
\[
\sum b_i \le \int \sup_j[b_j g_j(x)] \, d\mu\,.
\]
Then there are disjoint sets $E_j$ such that
\[
1\le \int_{E_j} g_j\, d\mu\,.
\]
\end{theorem}

Using this theorem, we see that \eqref{alCarlDual2} implies the existence of measurable subsets $E_Q\subset Q$ with  $Q\in \cD$, such that
\begin{equation}
\label{choice}
\al_Q\cdot \mu(Q)^2 \leq \mu(E_Q),\quad E_Q\cap E_{Q'} =\emptyset\,.
\end{equation}

We deduced \eqref{choice} from \eqref{alCarl1}. By now this has been done in various situations many times, see again \cite{Verb} (for cubes), \cite{TH} (for arbitrary Borel sets), and \cite{AB} for geometric proof of this geometric fact.  

The deduction of \eqref{choice} from \eqref{alCarl1} is called ``sparsity  property from Carleson property deduction".

Now we are ready to give the proof of Theorem \ref{Dor}.

\begin{proof} $(\Rightarrow)$ Suppose measure $\mu$ has no atoms and it is a bad measure. Then we can find a sequence $\al=\{\al_Q\}_{Q\in\cD}$ such that  condition \eqref{alC} is satisfied with constant $C=1$ (we can always normalize) but embedding \eqref{alEmb} is false. Then for every positive $K$ there exists a positive function $\psi$ such that 
\[
\sum_{Q\in \cD} \big(\int_Q \psi\, d\mu\big)^2 \al_Q\geq K\int \psi^2\, d\mu,
\]
or
\[
\sum_{Q\in \cD} \langle \psi\rangle_{Q, \mu}^2 \al_Q\cdot \mu(Q)^2\ge K\int \psi^2\, d\mu
\]

Let $\cF$ be the sub-collection of  $\cD$ such that $\al_Q\neq 0$ for $Q\in \cF$. We can use the sparsity condition \eqref{choice} to get

$$K\int \psi^2\, d\mu \leq \sum_{Q\in \cD} \langle \psi\rangle_{Q, \mu}^2 \al_Q\cdot \mu(Q)^2 \leq \sum_{Q\in \cF} \langle \psi\rangle_{Q, \mu}^2 \cdot \mu(E_Q)=$$

$$\sum_{Q\in \cF} \Big(\frac{1}{\mu(Q)}\int_Q\psi\hspace{0.1cm}{d}\mu\Big)^2  \cdot\mu(E_Q)\leq \sum_{Q\in \cF} \Big(\inf_{x\in Q}\cM_{\mu}\psi(x)\Big)^2\cdot \mu(E_Q)\leq$$
$$\sum_{Q\in \cF} \inf_{x\in Q}\Big(\cM_{\mu}\psi(x)\Big)^2\cdot \mu(E_Q)\leq\sum_{Q\in \cF} \inf_{x\in E_Q}\Big(\cM_{\mu}\psi(x)\Big)^2\cdot \mu(E_Q)\leq$$

$$\sum_{Q\in \cF} \int_{E_Q}\Big(\cM_{\mu}\psi(x)\Big)^2\hspace{0.1cm}{d}\mu\leq\int_{[0,1]^2}\Big(\cM_{\mu}\psi(x)\Big)^2{d}\mu$$

since the sets $E_Q$ are disjoint.

Hence
\begin{equation}
\label{MmuK}
\int \big(\cM_\mu\psi \big)^2\, d\mu \geq K\int \psi^2\, d\mu\,.
\end{equation}
Since $K$ can be arbitrary large the operator $\cM_\mu$ is not bounded in $L^2(\mu)$.

\bigskip

$(\Leftarrow)$
Now assume that operator $\cM_\mu$ is not bounded in $L^2(\mu)$. Then, for an arbitrary $K>0$ there exists a positive function $\psi$ such that \eqref{MmuK} holds. 

By approximating $\psi$ with an increasing sequence of positive, simple functions $\psi_n$ we can find $n_0\in \mathbb{N}$ such that \eqref{MmuK} holds with $\psi$ replaced by $\psi_{n_0}$. Also, we can find $N=N(n_0)\in\mathbb{N}$ so that the function $\psi_{n_0}$ is constant on each dyadic sub-square of $Q_0$ of size $2^{-N}\times 2^{-N}$. Given these, in the definition of $\cM_\mu\psi_{n_0}$ we replace $\sup$ by $\max$. For simplicity we omit writing the sub-script $n_0$.

Now we group the elements of $Q_0=[0,1]^2$ into the sets $A_Q$, $Q\in\cD$ as follows: $x\in A_Q$ if

$$\max_{\substack{R\in \cD \\ R \ni x}}\frac{1}{\mu(R)}\int_R \psi {d}\mu = \frac{1}{\mu(Q)}\int_Q \psi {d}\mu$$

Of course there may be more than one set $A_Q$ for each $x$. We would like to make these sets disjoint. Let ${Q_1}$, ${Q_2}$,.. be an enumeration of the dyadic rectangles $Q$, $Q\in \cD$. Then consider the sets $A_{Q_i}':=A_{Q_i}\setminus \bigcup\limits_{j=1}^{i}A_{Q_j}$. Obviously the sets $A_{Q_i}'$ are disjoint and their union is $Q_0=[0,1]^2$.

Hence, we have

$$\int_{Q_0}\Big(\cM_{\mu}\psi(x)\Big)^2{d}\mu=\int_{Q_0}\Bigg(\max_{\substack{R\in \cD \\ R \ni x}}\frac{1}{\mu(R)}\int_R \psi {d}\mu\Bigg)^2{d}\mu=\sum_{i}\int_{A_{Q_i}'}\Big(\frac{1}{\mu(Q_i)}\int_{Q_i} \psi {d}\mu\Big)^2 {d}\mu =$$

$$\sum_{i}\frac{\mu(A_{Q_i}')}{\mu(Q_i)^2}\Big(\int_{Q_i} \psi \hspace{0.1cm}{d}\mu\Big)^2:=\sum_{i}\al_{Q_i}\Big(\int_{Q_i} \psi \hspace{0.1cm}{d}\mu\Big)^2$$

In other words, we constructed a sequence $\al=\{\al_{Q_i}\}$, $Q_i\in \cD$ for which the embedding \eqref{alEmb} does hold with a very large constant, but the Carleson condition holds with constant $1$. But, this sequence, as well as the function $\psi$ depend on $K$. To get a general sequence and function we do start with an $i\in\mathbb{N}$ and take $K=4^i$. Lets omit enumerating the dyadic rectangles in $[0,1]^2$ for now. Given the above,  $\exists \psi_i$ and $\{\al_{Q}^i\}_{Q\in\cD}$ such that 

$$4^i \int\limits_{[0,1]^2}\psi_i^2{d}\mu \hspace{0.1cm}\leq \hspace{0.1cm} \sum\limits_{Q\in \cD}\al_{Q}^i\Big(\int_{Q} \psi_i \hspace{0.1cm}{d}\mu\Big)^2 $$ and re-normalizing we get 

\begin{equation}\label{iemb}
2^i  \hspace{0.1cm}\leq \hspace{0.1cm} \sum\limits_{Q\in \cD}\al_{Q}^i\Big(\int_{Q} \phi_i \hspace{0.1cm}{d}\mu\Big)^2 
\end{equation}

for  $\phi_i=\dfrac{\psi_i}{2^{i/2}\|\psi_i\|_{L^2([0,1]^2)}}$.

Now let $\al_Q:=\sum\limits_{i=1}^{\infty}2^{-i}\al_Q^i$ and $\phi=\sum\limits_{i=1}^{\infty}\phi_i$. Notice that $\phi\in L^2([0,1]^2)$ and the sequence $\al=\{\al_Q\}_{Q\in \cD}$ satisfies the Carleson condition: Let any sub-collection $\cS$ of $\cD$. There is an indexing set $J_{\cS}=\{j_1,j_2,...\}$ which enumerates the dyadic rectangles in $\cS$. By the disjointness of $A_{Q_j}'^{i}$'s (for each fixed $i$) we have

$$\sum_{Q\in \cS}\al_{Q} \mu(Q)^2=\sum_{j\in J_{\cS}}\al_{Q_j} \mu(Q_j)^2=\sum_{j\in J_{\cS}} \sum\limits_{i=1}^{\infty}2^{-i}\al_{Q_j}^i \mu(Q_j)^2 = \sum\limits_{i=1}^{\infty}2^{-i} \sum_{j\in J_{\cS}} \mu(A_{Q_j}'^{i})= $$

$$\sum\limits_{i=1}^{\infty}2^{-i}\mu(\bigcup\limits_{j\in J_{\cS}} A_{Q_j}'^{i}) \leq \sum\limits_{i=1}^{\infty}2^{-i} \mu(\bigcup\limits_{Q \in {\cS}} {Q})=\mu(\bigcup\limits_{Q \in {\cS}} {Q})  $$

Finally, using \eqref{iemb} we get

$$\sum\limits_{Q\in\cD}\al_Q \Big(\int_{Q} \phi \hspace{0.1cm}{d}\mu\Big)^2=\sum\limits_{Q\in\cD}\sum\limits_{i=1}^{\infty}2^{-i}\al_Q^i\Big(\int_{Q} \phi \hspace{0.1cm}{d}\mu\Big)^2 \ge \sum\limits_{i=1}^{\infty}2^{-i}\sum\limits_{Q\in\cD}\al_Q^i\Big(\int_{Q} \phi_i  \hspace{0.1cm}{d}\mu\Big)^2 = \infty$$

\end{proof}

\end{document}